%Paper: alg-geom/9501011
%From: Rikard B|gvad <rikard@matematik.su.se>
%Date: Fri, 20 Jan 1995 10:09:08 +0100

\input amstex
%A4
\hsize 124.2mm
\vsize 196.2mm
\baselineskip=13pt

\magnification =1200

%MAKRON

\def\dlyft1#1{{\Cal D}^{(1)}_{\tilde#1}}
\def\surj{\to\kern-.75em \to }

\def \o{o \kern -.9 em \slash}
\def \O#1{{\Cal O}_{#1}}

\def \o#1{{\Cal O}(#1)}
\def \L{{\Cal L}}

\def \R {\operatorname{\bold R}}

\def \spec {\operatorname {spec}}
\def \Pic #1 {\operatorname {Pic} #1}

\input amsppt.sty
\topmatter

\title Some homogeneous coordinate rings that are
Koszul algebras \endtitle

 \author Rikard B\"ogvad\endauthor
\abstract {Using reduction to positive characteristic
and the method of Frobenius splitting of diagonals, due to Mehta and
Ramanathan, it is shown that homogeneous coordinate rings for either
proper and smooth toric varieties or Schubert varieties are Koszul algebras.
}\endabstract

\endtopmatter

\vskip 2cm
\subheading {1.Introduction}
All varieties will be assumed to be defined and proper over an algebraically
closed field $k$. Assume that $\L$ is a line bundle on an
algebraic variety $X$,
and define the graded $k$-algebra $$R(\L):=\oplus_{n\geq 0}\Gamma (X,
\L^{\otimes n}.)$$ Since $X$ is proper this algebra is finite dimensional in
each degree as a $k$-vector space,  and it is the integral closure of the
homogeneous coordinate ring of the image of $X$ under the map defined by
$\L $. The striking result of  Mumford that $R(\L^{\otimes n})$ has quadratic
relations if $n$ is large enough has been generalized in several ways, for
example resulting in an analysis of the degrees of the syzygies of $R(\L)$
considered as a quotient of a minimal polynomial algebra (cf. [EL]). Another
generalization is the theorem of J. Backelin [B] that for $n$  large enough
$R(\L^{\otimes n})$ is a Koszul algebra (See Def 2.1.). G.R.Kempf has noted in
[K2] that for certain types of varieties (e.g. abelian varieties, Grassman
varieties, and curves) $R(\L)$ is Koszul except for very few $\L$, and he
suggests the problem of finding more examples .  Our first result is the
following.

\proclaim {Theorem 1}  Let $X$ be
a nonsingular proper toric variety, and suppose that $\L$ is an ample line
bundle on $X$. Then $R(\L)$ is a Koszul algebra. \endproclaim

In the case of homogeneous spaces we obtain the slightly stronger result

\proclaim {Theorem 2} Let $X$ be
a Schubert variety in a homogeneous space $G/P$ of a reductive algebraic
group
 modulo
a parabolic subgroup, and suppose that $\L$ is an effective line
bundle on $G/P$. Then $R(\L)$ is a Koszul algebra. \endproclaim

It is known (cf.[B,R]) that $R(\L)$ is generated
by $\Gamma\L$ in either case so $R(\L)$ is the homogeneous coordinate ring of
$X$ in the embedding given by $\L$. A Koszul algebra always has
quadratic relations, so this is a generalization of the results of Ramanathan
[R]  -- homogeneous spaces -- and [B] -- toric varieties. In [K], it is shown
that algebras with a straightening law whose discrete algebra is defined by
quadratic monomials are Koszul algebras. These algebras include, as stated
there, at least embeddings of Grassmannians by Pl\"ucker coordinates, but I
ignore whether they include all homogeneous varieties of the theorem.

There is also a relative version of Theorem 2, simplifying proofs, and
generalizing the statement that Schubert varieties are linearly defined in an
embedding of the corresponding Borel variety given by an effective line bundle
$\L$ (cf.[R],[KR]). A ring homomorphism is said to be 1-linear if it is is
relatively Koszul (see Def.2.1.).

\proclaim {Theorem 3} Let $X$ be as in Theorem 2. Suppose
that $Y\subset X$ is a
Schubert variety and that $\L$ is an effective line
bundle on $X$. Let $S:= R(Y,\L)$ be the restriction to
$Y$ of $R=R(X,\L)$. Then the surjective map $s:R\to S$ is
1-linear (cf. Def 2.1).\endproclaim

The above results follow by
generalizing the methods used by Ramanathan to prove that homogeneous
coordinate rings of Schubert varieties have quadratic relations to higher
syzygies. This generalization is contained in the following
theorem. Assume that
$Y\subset X$ is a closed subvariety, and  define the {\it partial diagonal}
$D^s_{i,i+1}\subset Y\times X\times \ldots\times X=Y\times X(s)$ ($s$ factors
$X$), as the set $$\{(x_0,x_1,\ldots,x_s),\ x_0\in Y,\
x_k\in X,\ 1\leq k\leq s,\
x_i=x_{i+1}.\}.$$

\proclaim {Theorem 4} Assume that the characteristic of $k$
is positive and that
$Y\subset X$ is a closed subvariety of the proper variety
$X$, possibly the empty
set (this corresponds to the Koszul property, see section 2). Assume that each
$Y\times X(n)$, $n=1,2,\ldots.$ is Frobenius split and that all diagonals
$\Delta_{ii+1}^n$ for $
i=0,\ldots,n-1$, are simultaneously compatibly split as closed subvarieties of
$X(n)$. Let $S:= R(Y,\L)$ be the restriction to $Y$ of
$R:=R(X,\L)$. Then the map
$s:R\to S$ is 1-linear. \endproclaim

\medskip

 The idea of the proof of Theorem 4 is to find conditions for
the existence of a (partial) homotopy of the identity on a certain bar complex
with the zero map,  showing that the homology is concentrated in the right
degrees.  These conditions, infinite in number, are
formulated as the acyclicity
in the right dimension of a kind of model, as is done in,
for example, the proof
of the Eilenberg-Zilber theorem in elementary algebraic topology [M,
VIII,7-8]. (See Theorem 5.)  Since we are working with the bar complex they are
then easily rephrased in terms of restriction maps of
cohomology of line bundles
to  partial diagonals. Then the results of Mehta and Ramanathan on compatible
splitting show that the resulting conditions are satisfied if the partial
diagonals are compatibly split.

Since our conditions for having a Koszul algebra unfortunately are infinite in
number, Backelin's theorem does not follow from them
(examples have been constructed by J.-E Roos [Ro] showing that the right
behaviour for an arbitrary graded ring of a finite number of torsion groups
does not suffice to prove that it is a Koszul algebra, even if the number of
generators of the algebra is bounded). Thus one might hope that a closer
analysis of a resolution of a homogeneous algebra would give much sharper
results, valid generally.
In a sequel a generalization of the above results to multi-cones as in [KR ]
will be given.

I have been told that results similar to those in the
present note on homogeneous
varieties have been obtained by M.S.Ravi.

I would like to thank J.Backelin, T.Ekedahl, R.Fr\"oberg, C.L\"ofwall and J.-E
Roos for helpful conversations on these topics.

\subheading {2.1 The bar complex} Let $k$ be a field (of
arbitrary characteristic). We will only
consider finitely generated commutative graded $k$-algebras
 $$
R=\oplus_{i\geq 0} R_i
$$
who are {\it connected}, i. e. $R_0=k$. Homomorphisms
between such graded rings will always be
homogeneous of degree 0. Let $R^+:=\oplus_{i> 0} R_i$ be the irrelevant ideal.
 The normalized standard complex (often
called the {\it bar complex}) of $R$ has underlying vector space
$$
N.=\bigoplus_{n\geq 0}R\otimes_k
(R^+)^n.
$$
Here $(R^+)^n$ denotes the tensor product over $k$ of $R^+$ with itself $n$
times (letting $(R^+)^0=k$) (Cf. [CE, X.2]). A typical
element of the complex is denoted
by $$\lambda_0[\lambda_1,\lambda_2,\ldots,\lambda_n ]:=
 \lambda_0\otimes\lambda_1\otimes\lambda_2\otimes\ldots\otimes\lambda_n
\in N_n.$$
The $R$-linear differential $d_n: N_n\to N_{n-1}$ of
the complex is given by
$$
\split d_n[\lambda_1,\lambda_2,\ldots,\lambda_n
]:&=\lambda_1[\lambda_2,\ldots,\lambda_n]+\sum
_{0<i<n}(-1)^{i}[\lambda_1,\lambda_2,\ldots,
\lambda_i\lambda_{i+1},\ldots.\lambda_n]\\
\endsplit\tag1
$$
 With this differential the complex is an $R$-free
resolution of $k$. Observe that
it is bigraded; in addition to its grading as a complex ---
its {\it homological degree} --- it has
a second
grading induced from the grading of $R$, called the {\it internal
degree}. The  homogeneous piece of $N.$  of bidegree $(i,j)$
will be denoted by $N_{ij}$.  The
differential has bidegree $(-1,0)$ as a homogeneous map in
this bigrading. If $M$ is a graded
$R$-module, the induced differential of $M\otimes_RN.$ is
still homogeneous, and hence the bigrading
passes to the homology $$H(M\otimes_RN.)=Tor_.^R(M,k)
=\bigoplus_{i,j\geq 0}Tor_{ij}^R(M,k).
$$

 We now introduce more notation. Assume that
S$:=\oplus_{i\geq 0} S_i$, is a  graded
connected $R$-algebra. Note that it is immediately clear
from the construction of the
bar complex that $(S\otimes_RN)_{ij}=(S\otimes_k
(R^+)^i)_j=0$, if $j<i$ and that hence in this case also
$Tor_{ij}^R(S,k)=0$.  The {\it diagonal
part} $D(S):=\bigoplus_{n\geq 0}S_0\otimes_k (R_1)^n,$ is the subspace of
$S\otimes_RN.$ consisting in each homological degree of the
part of minimal internal degree. The graded complement
$V(S):=\bigoplus_{n=0}S^+\otimes_k
(R^+)^n,$ will be called the {\it off-diagonal part} of the complex
$S\otimes_RN.$ It is a subcomplex and contains all
boundaries, by reason of degree. Hence there is
an injection $$Z(S\otimes_RN.) \cap D(S)  \cong \bigoplus_{i\geq
0 }Tor^R_{ii}(S,k)\hookrightarrow H(S\otimes_RN.).\tag 2$$

Denote the free
$R$-module of rank $d$ with generators in degree $j$ by $R[j]^{d}$.
 The dimensions
of the torsion groups $d_{ij}=dim_kTor_{ij}^R(S,k)$ may also be interpreted in
terms of ranks of
the minimal free graded resolution $F_.$ of $S$. Namely it is true that
$F_i=\oplus_jR[j]^{d_{ij}}$. This motivates the first part of the following
(fairly standard) terminology:
\proclaim {Definition 2.1}(Cf. [BF], [EG]). If $S=:R/I$ is a
quotient of $R$, such
that $Tor_{ij}^R(S,k)=0$ if $i\neq j$, then $S$ is said to have a 1-linear
$R$-resolution, or the map $R\to S$ is said to be 1-linear. If $k=R/R^+$ has a
1-linear resolution then $R$ is called a Koszul algebra.
\endproclaim
Using the previous terminology it is clear that $R\to S$ is
1-linear iff the inclusion
$V(S)\subset S\otimes_RN.$ induces the zero map on homology.
In the more enthusiastic terminology of Kempf [K] one says
that $R$ is wonderful if
$R$ is a Koszul algebra, and that $S$ is an awesome module if it has a 1-linear
resolution. Note that a Koszul algebra $R $ is generated by
$R_1$ (since $Tor^R_{1j}(k,k)=0$ if
$j\neq 1$, see e.g. [Le]) and if it is presented as a graded quotient of the
polynomial algebra $k[R_1]$, it will have relations in
degree 2 (since $Tor^R_{2j}(k,k)=0$ if $j\neq
2$, loc.cit.). It is also clear that a ring $S$ with an
1-linear resolution will be a quotient of $R$
by an ideal generated by elements of degree 1.

We have use for another concept, due to G.Levin [L], from
the theory of the homology of
commutative rings.

\proclaim {Definition 2.2} A homomorphism  $\eta: R\to S$ that
induces a surjection $Tor^R_.(k,k)\to Tor^S_.(k,k)$ is
called a large homomorphism.\endproclaim
Some examples are given by the following lemma.
\proclaim {Lemma 2.3} a) A retraction
$\eta:R\to S$  onto a graded subalgebra $S\subset R$ is a large homomorphism.

b) Suppose that $R \to S=R/I$ is 1-linear. Then
it is large.

 c) If $R $ is a Koszul algebra and $\eta:R\to S$ is a large map, then  $
S$ is also a Koszul algebra, and furthermore $\eta:R\to S$ is then 1-linear.
 \endproclaim
Part a) is wellknown [L, Theorem 2.3] and immediate (see the
similar proof in Lemma 3.1.2 in [B].)
 The remaining statements are a generalization of [Ba, Lemma 2.3]
(Backelin assumes that both rings involved are Koszul
algebras). The proof uses the
following characterization of large homomorphisms $\eta:R\to S$, given in [L,
Theorem 1.1]. Namely, $\eta $
is large iff the map  $$p_*:Tor^R(S,k)\to Tor^R(k,k), $$ induced by the natural
surjection $p:S\to S/S^+=k$ is an injection (Levin studies local rings but his
arguments are true in our graded situation). Let us prove that
this criterion is true if $\eta:R\to S$ is 1-linear, i.e.
that $Tor^R_{ij}(k,k)=0$ if $i\neq
j$. Consider the bar complex $N_.$ of $R$. Then $p$ induces $p_*$:
$$
S\otimes_R N_.\to k\otimes_R N_.
$$
This map has bidegree $(0,0)$, so it maps in particular the
diagonal part $$D(S):=\bigoplus_{n\geq
0}S_0\otimes_k (R^+)^n,$$ isomorphically to
$D(k):=\bigoplus_{n\geq 0}k\otimes_k
(R^+)^n,$ and since, as mentioned above, by reason of
degree, $D(k)$ contains no boundaries it
follows that $p_*$ is injective on the diagonal part of the homology:
$$
p_*:Z(S\otimes_R N_.)\cap D(S)  \cong \bigoplus_{i\geq
0 }Tor^R_{ii}(S,k)\hookrightarrow Z(k\otimes_R N_.)\cap
D(k)  \cong \bigoplus_{i\geq
0 }Tor^R_{ii}(k,k).
$$
(Note that this known result is valid for arbitrary rings).
If $R\to S$ is 1-linear
the only non-zero homology is the one on the diagonal, and
thus the map is large.
Finally the last part of the lemma is an immediate
consequence of the definitions, and the
previously used characterization of large homomorphisms.

\subheading {2.2 The bar complex of a homogeneous coordinate ring}
The following situation will be our main object of study.
Suppose that $\L$ is a line bundle on
the proper variety $X$. Form, as in the introduction the
graded connected algebra
$$
R=R(X,\L)=k\oplus R^+.
$$
Let $Y$ be a subvariety of $X$, restrict $\L$ to $Y$ and
consider in the same way the $k$-algebra
$S:=S(\L)$. There is a canonical restriction homomorphism $R\to S$.
 Consider the bar complex  $N.$ of $R$ tensored with $S$.
$$
S\otimes_kN.=\bigoplus_{n\geq 0}
S\otimes_k(R^+)^n\cong\bigoplus _{1\leq n,0\leq i_0,1\leq
i_1,\ldots i_n}\Gamma(Y\times
X(n),\L^{i_0}\times \ldots \times \L^{i_n}),
$$
Here $Y\times X(n)$ denotes the product of $Y$ with
with $n$ copies of $X$, and we have used the canonical identification  $$
\Gamma(Y\times
X(n),\L^{i_0}\times
\ldots \times \L^{i_n})\cong
\Gamma(Y,\L^{\otimes i_0})\otimes_k\ldots
\otimes_k\Gamma(X,\L^{\otimes i_n}).
$$

The differential $d_n: N_n\to N_{n-1}$ of the complex may
be given as the sum
$$
\split d_n[\lambda_1,\lambda_2,\ldots,\lambda_n
]:&=\lambda_1[\lambda_2,\ldots,\lambda_n]+\sum
_{0<i<n}(-1)^{i}[\lambda_1,\lambda_2,\ldots,
\lambda_i\lambda_{i+1},\ldots.\lambda_n]\\
&=d^n_{01}
[\lambda_1,\lambda_2,\ldots,\lambda_n]+\sum_{0<i<n}(-1)^id^n_{ii+1}
[\lambda_1,\lambda_2,\ldots,\lambda_n].\endsplit\tag1
$$
 of restriction maps $$\split d^n_{ii+1}:&S\otimes_k(R^+)^n \supset
\Gamma(Y\times X(n),\L^{i_0}\times \ldots \times \L^{i_n})\to
\Gamma(\Delta^n_{ii+1},\L^{i_0}\times \ldots \times \L^{i_n})\\
&\cong \oplus
\Gamma(Y\times X(n-1),\L^{i_0}\times \ldots \times
\L^{i_0}\otimes \L^{i_{i+1}} \times
\L^{i_n})\subset  S\otimes_k(R^+)^{n-1}\endsplit$$ for
$i=0,\ldots n-1,$ on cohomology from
$Y\times X(n)$ to partial diagonals $\Delta^n_{ii+1}$
consisting of those points in $Y\times X(n)$
that have coordinate $i$ equal to coordinate $i+1$.
 These maps will be called {\it partial
differentials}.

With joyfully Bourbakistic nostalgia we extend the
definition of homogeneous coordinate rings to the empty
subset, and defining $R(\emptyset
,\L):=k$ if $\L$ is the restriction of a line bundle $\L$ on
$X$, we obtain the quotient $k$ as the
restriction to a closed subvariety. This gives us the added
comfort to be able to
treat uniformly possibly 1-linear restrictions of $R$ as
well as the question whether $R$ itself is a
Koszul algebra.

It will be convenient to study these maps in a more general
context. Let ${\Cal I}^n_{ii+1}$ be the
ideal of $$\Delta^n_{ii+1}\subset Y\times X(n),$$ and denote
the inclusion by $i$. Assume that
$\L_\alpha:= \L_{0}\times\ldots \times \L_{n},$ is an
invertible sheaf on $Y\times X(n)$, which
is the product of the invertible sheaves $\L_i$. There is then a
short  sequence.
$$ \split\Gamma(Y\times X(n),{\Cal I}^n_{ii+1}\L_\alpha)\hookrightarrow&
\Gamma(Y\times X(n),\L_\alpha)\to\\
 \to\Gamma(\Delta^n_{ii+1},i^*\L_\alpha)&\cong
\Gamma(Y\times X(n-1),\L_{0}\times\ldots \times \L_i\otimes
\L_{i+1} \times\ldots\times
\L_{n}).\endsplit$$
If $\L_j=\L^{i_j}$, the map composed of the last two
maps in the short sequence is precisely the differential
$d^n_{ii+1}$, in the complex $S\otimes_RN.$ described above.
For this reason, we will
continue to call the composite map $d^n_{ii+1}$ in the
general case, suppressing
any reference to the line bundles are involved. The first part of the following
lemma is clear. The second part is a consequence of the fact that taking global
sections commutes with finite inverse limits.
 \proclaim
{Lemma 2.2} If $X$ is a proper variety and $\L_\alpha$ and  $d^n_{ii+1}$ are as
above then $$Kerd^n_{ii+1}= \Gamma(Y\times X(n),{\Cal I}^n_{ii+1}\L_\alpha).$$
Furthermore $$\cap_{i\in I}Kerd^n_{ii+1}=\Gamma(Y\times X(n),
\cap_{i\in I}{\Cal I}^n_{ii+1}\L_\alpha),\ {\text if} \
I\subset\{0,\ldots,n\}. $$\endproclaim

\subheading{3. Models for constructing a homotopy}
Let $\L_\alpha=
\L_{0}\times\ldots \times \L_{n+1}$, for $n\geq 0$, be a line bundle on
$Y\times X(n+1)$. We will now construct an associated {\it model complex}
$K=K(n,\L_\alpha).$ It is  zero in all (homological) degrees except possibly
$n-1,n,n+1 $ :
 $$
K=K_{n-1}\oplus K_n\oplus K_{n+1}.$$
and
 $$ \align
 &K_{n+1}:= \Gamma(Y\times X(n+1),\L_\alpha)\\
&K_n:=\bigoplus_{0\leq i\leq n} \Gamma(Y\times X(n),\L_{0}\times\ldots
\times \L_i\otimes \L_{i+1} \times\ldots\times
\L_{n+1})\\
 &K_{n-1}:=\bigoplus_{0\leq k<i\leq
n} \Gamma(Y\times X(n-1),\L_{0}\times\ldots\times
\L_k\otimes \L_{k+1}\times \ldots \times
\L_i\otimes \L_{i+1} \times\ldots\times \L_{n+1})
\endalign$$
The differential restricted to $K_{n-1}$ is zero, while it
is the same as in the
standard complex acting on $K_n$ or $K_{n+1}$.  Thus, if a typical element of
$K_{n+1}$ or $K_n$ is denoted by the symbol
$[\lambda_0,\ldots,\lambda_s],$ for $\
s=n-1,$ or $n,$ where each  $\lambda_i\in\Gamma(X,M)$ for some $M=\L_j$ or
$M=\L_j\otimes \L_{j+1}$, then  formula 2.(1) precisely describes the
differential. (Using the partial differentials described as
a preliminary of lemma
2.1). For the definition to be intelligible for $n=0$ we consider $X(0)$ to be
$\spec k$, and $X(-1):=\emptyset$. Then $K(0,\L_0\times \L_1)$ is the complex
$$
\Gamma(Y\times X,\L_0\times \L_1)\cong\Gamma(Y,\L_0)\otimes\Gamma(X,
\L_1) \to \Gamma(Y,\L_0\otimes \L_1)\to 0,\tag 1 $$
and the only non-zero differential $d_1=d_1^{01}$
corresponds to multiplication of
global sections.

We will now inductively describe the homology of the models.
\proclaim {Lemma 3.1} Assume that
$\L_\alpha=\L_{0}\times\ldots \times \L_{n+1}$ is an invertible sheaf on
$Y\times X(n+1)$ and let $\tilde \L_\alpha$ be the
 sheaf $\L_{0}\times\ldots \times \L_{n} $ on $Y\times X(n)$. The partial
differential

 $$d_{n+1}^{nn+1}:\  \
A:=\Gamma(Y\times X(n+1),\L_\alpha)\to B:=\Gamma(Y\times
X(n),\L_{0}\times\ldots \times \L_n\otimes\L_{n+1})$$
induces a
map $\delta$ of subspaces:
$$
 \delta: A\supset \bigcap_{0\leq i\leq n-1}\operatorname {ker}d_{n+1}^{ii+1}\to
\bigcap_{0\leq i\leq n-1}\operatorname {ker}d_{n}^{ii+1}\subset B,
$$
and there is an exact sequence
$$
 cok\delta\hookrightarrow H_n(K(n,\L_\alpha))\to
H_{n-1}(K(n-1,\tilde \L_\alpha))\otimes
\Gamma(X,\L_{n+1}).$$If $n=0$ there is an
isomorphism
$$
H_0(K(0,\L_\alpha)\cong cok \delta=cok d^1_{01}.
$$ \endproclaim

\demo {Proof} The case $n=0$ is immediately clear from the description (1) of
$H_0(K(0,\L_0\times \L_1)$.

 Note that, forgetting differentials, $K(n,\L_\alpha)$ is the
direct sum of vectors paces of global sections of the form
$\Gamma(Y\times X(s),M_0\times\ldots\times M_s)$, where
$M_i,\ 0\leq i\leq s$, are
invertible sheaves on $Y$ or $X$ and $s=n-1,n$ or $n+1$. Let
$V$ and $W$ be the sub
vectors paces of $K(n,\L_\alpha)$ consisting of the direct sums of those
$\Gamma(Y\times X(s),M_1\times\ldots\times M_s)$ with
$M_s=\L_{n+1}$, respectively
$M_s\neq \L_{n+1}$. Then $K(n,\L_\alpha)=V\oplus W$, and $W$ is clearly a
subcomplex.
It is clear that there are isomorphisms of vectors paces
$$V\cong K(n,\L_\alpha)/W\cong K(n-1,\tilde
\L_\alpha)\otimes\Gamma(X,\L_{n+1}).\tag2$$ The last vectors pace may be
considered as a complex, namely as the tensor product (over
$k$) of $K(n-1,\tilde
\L_\alpha)$ with the complex $\Gamma(X,\L_{n+1})[-1]$,  which has differential
zero and is concentrated in degree 1. (Here, as usual, $M.[-1]$ of a complex
$(M.,d)$ denotes the complex, which in (homological) degree $m $ has
$M[-1]_m=M_{m-1}$ and differential $-d$.) Then the final vector space
isomorphism of (2) is in fact an isomorphism of complexes.
Now take the long homology sequence belonging to the short exact sequence of
complexes

$$W\hookrightarrow K(n,\L_\alpha)\overset\beta\to\surj K(n-1,\tilde
\L_\alpha)\otimes\Gamma(X,\L_{n+1})[-1].$$
We get the exact sequence

$$\split H_{n+1}(K(n-1,\tilde \L_\alpha)\otimes\Gamma(X,\L_{n+1})[-1]&\to
H_n(W)\to\\ \to H_n(K(n,\L_\alpha))\to &H_n(K(n-1,\tilde
\L_\alpha)\otimes\Gamma(X,\L_{n+1})[-1]),\endsplit\tag 3$$
Since $W_{n+1}=0$ and
$W_n=\Gamma (Y\times X(n),\L_0\times\ldots \times \L_{n-1}\times \L_{n}\otimes
\L_{n+1})$, it is clear that
$$
H_n(W)=\operatorname {kerd}_n\cap W_n=ker(\Sigma_{i=0}^{n-1}
\operatorname {kerd}_n^{ii+1})\cap W_n=\bigcap_{0^\leq i\leq
n-1}\operatorname {kerd}_n^{ii+1}\cap W_n. $$
Similarly

 $$ H_{n+1}(K(n-1,\tilde
\L_\alpha)\otimes\Gamma(X,\L_{n+1})[-1])= 
\bigcap_{0\leq i\leq n-1}\operatorname {kerd}_{n+1}^{ii+1},\tag4$$ (the
intersection taken in $A$ of the lemma),  and
 $$H_{n}(K(n-1,\tilde
\L_\alpha)\otimes\Gamma(X,\L_{n+1})[-1])=H_{n-1}(K(n-1,\tilde
\L_\alpha))\otimes\Gamma(X,\L_{n+1}),
$$
so it
suffices to check that the connecting map in the long exact
sequence (3) is given
by $d_{n+1}^{nn+1}$. But this follows, from (4), since,
using the canonical lifting
$K(n,\L_\alpha)_{n+1}/W_{n+1}\cong V_{n+1}\subset K(n,\L_\alpha)_{n+1}$, the
connecting map is $d_{n+1}$ restricted to $V$, and since the maps
$d_{n+1}=\Sigma_{i=0}^n d_{n+1}^{ii+1}$ and $d_{n+1}^{nn+1}$
evidently coincide on
$\cap_{0\leq i\leq n-1}\operatorname {kerd}_{n+1}^{ii+1}$.

\enddemo

  By induction it is now possible to give reasonable geometric
conditions for the models to have middle degree homology zero.

\proclaim {Lemma 3.2} Assume that the following conditions are satisfied
by an invertible sheaf $\L_\alpha$ on $Y\times X(n+1):$

i)  $H^1(Y\times X(n+1),\L_\alpha)=0$.

ii) $\Gamma(Y\times X(n+1),\L_\alpha)\surj
\Gamma (\bigcup_{i=0}^{n-1}\Delta^{n+1}_{ii+1},\L_\alpha)$

iii) $\Gamma(Y\times X(n+1),\L_\alpha)\surj
\Gamma (\Delta^{n+1}_{nn+1},\L_\alpha)$

Then $cok\delta$ (see the preceding lemma) is zero. If the conditions i) to
iii) are satisfied for all $\L_0\times
\ldots \times \L_{m+1}$ where $0\leq m\leq n$ then $H_n(K(n,\L_\alpha))=0$

\endproclaim

\demo {Proof}Recall that ${\Cal I}^{n+1}_{ii+1}$ is the
ideal in $Y\times X(n+1)$
of the partial diagonal $\Delta^{n+1}_{ii+1}$ (and that $n+1$ is here not a
power but an index).
 Observe that conditions i) and
ii) of the lemma, together with the exact sequence
 $$ \split\Gamma(Y\times X(n+1),\L_\alpha)\surj
&\Gamma (\bigcup_{i=0}^{n-1}\Delta^{n+1}_{ii+1},\L_\alpha)\to\\ \to
&H^1(Y\times X(n+1),\bigcap_{i=0}^{n-1}{\Cal I}^{n+1}_{ii+1}\L_\alpha)\to
H^1(Y\times X(n+1), \L_\alpha)\endsplit $$
imply that $H^1(Y\times X(n+1),\bigcap_{i=0}^{n-1}{\Cal
I}^{n+1}_{ii+1}\L_\alpha)=0$.
(Note that $\bigcap_{i=0}^{n-1}{\Cal I}^{n+1}_{ii+1}$ is the ideal of
$\bigcup_{i=0}^{n-1}\Delta^{n+1}_{ii+1}$).

There is a short exact sequence
$$ \bigcap_{i=0}^{n}{\Cal
I}^{n+1}_{ii+1}\hookrightarrow \bigcap_{i=0}^{n-1}{\Cal
I}^{n+1}_{ii+1}\oplus{\Cal
I}^{n+1}_{nn+1}\surj {\Cal J}:=\bigcap_{i=0}^{n-1}{\Cal I}^{n+1}_{ii+1}+ {\Cal
I}^{n+1}_{nn+1}\subset  \O
{Y\times X(n+1)}.\tag5$$
Taking the tensor product of (5) with $\L_\alpha $ (over $\O{}  =\O
{Y\times X(n+1)}$) and taking the long cohomology sequence shows that
$$
\Gamma(Y\times X(n+1),(\bigcap_{i=0}^{n-1}{\Cal
I}^{n+1}_{ii+1}\oplus{\Cal I}^{n+1}_{nn+1})\L_\alpha)\surj \Gamma(Y\times
X(n+1),{\Cal J}\L_\alpha)\tag6
$$
is a surjection.
 Let ${\Cal I}_1:=\bigcap_{i=0}^{n-1}{\Cal I}^{n+1}_{ii+1}$, ${\Cal
I}_2:={\Cal I}^{n+1}_{nn+1}$ and $S:=Y\times
X(n+1)$, and consider the commutative diagram:
$$\matrix
\Gamma(S,{\Cal I}_1\L_\alpha\oplus{\Cal
I}_2\L_\alpha)&&&&\\ \downarrow&&&&\\
\Gamma(S,J\L_\alpha)&\hookrightarrow&\Gamma(S,\L_\alpha)&\to& \Gamma(S,\O{}
/J\L_\alpha)\\ \downarrow \epsilon&&\downarrow\eta&{\downarrow\ =}\\
\Gamma(S,J/{\Cal I}_2 \L_\alpha)&\hookrightarrow&\Gamma(S,\O{}
 /{\Cal I}_2\L_\alpha)&\to& \Gamma(S,\O{}  /J\L_\alpha)\\
{\downarrow \cong}&&{\downarrow \cong}&&\\
\Gamma(Y\times X(n),\bigcap_{i=0}^{n-1}{\Cal I}^{n-1}_{ii+1}
\L_\alpha)&\hookrightarrow&\Gamma(Y\times X(n),\L_0\times\ldots
\L_n\otimes \L_{n+1})&&\\\endmatrix $$
The map composed of the middle vertical maps is precisely the partial
differential $d_{n+1}^{nn+1}$. The composition of the
left hand vertical morphisms kills $\Gamma(Y\times
X(n+1),{\Cal I}_2\L_\alpha)$,and then
by lemma 2.2, the cokernel of this composed map equals  $cok\delta$. In view of
(5), this implies that $cok\delta\cong cok\epsilon$. By a diagram chase,
$cok\epsilon=0$ if $cok\eta=0$, since the right hand vertical map is injective.
This, however, is condition iii) of the lemma, which hence is proved.(The
last statement follows from the preceding lemma by an obvious induction).
\enddemo

\subheading{4.1. A criterion for a restriction map of homogeneous coordinate
rings to be 1-linear}

We will now study the homology of homogeneous coordinate rings, using the model
complexes.
The object of this paragraph is the following result.
 \proclaim{Theorem 5}Let $X$ be a proper variety, $\L$ a line bundle
on $X$ and $R:=R(X,\L)$. Suppose further that $S$ is the
restriction of $R$ to the
closed subvariety $Y$ of $X$, and that the conditions i-iii) of Lemma 3.2 are
fulfilled for all line bundles of the type $\L^\alpha:=\L^{\otimes
a_0}\times\ldots\times \L^{\otimes a_{n+1}}$ where $a_0\geq 0,$ and $\ a_i>  0$
if $i>0$. Then $R\to S$ is 1-linear. (In particular note that putting
$Y=\emptyset $ in Lemma 3.2. gives conditions for $R$ to
be Koszul.) \endproclaim

The proof will take the rest of the section. By Lemma 3.2
all higher homology is
zero for all model complexes:
$$
H_nK(n+1,\L^{\otimes a_0}\times\ldots\times
\L^{\otimes a_{n+1}})=0, \text{
if} \ n\geq1,\tag1
$$
 and $a_0\geq 0,\ a_i>  0$ if $i>0$.
 Let
$D=D(S)$ be the diagonal vector subspace of $K.:=S\otimes_RN.$ which contains
everything on the diagonal,  $$ D=\bigoplus
_{n\geq 0}\Gamma(Y\times
X(n),\O Y\times \L\times\ldots\times
\L).$$  Consider $D$ as a complex with the zero
differential. There is a surjective map of complexes $K.\to D$ (with kernel
the off-diagonal subcomplex $V:=V(S)$) (see section 2.1).
The theorem asserts that this map gives an injection after taking homology. To
prove this it suffices to prove that the injection $V.\subset K.$ is
homotopic to 0.

\subheading{4.2 Homotopy}

 First we will define a function which will give the target of
the homotopy.
 The function $\sigma$ is a map from sequences $(\alpha):=(a^0,\ldots,a^n)$ of
length $n+1$
 consisting of nonnegative integers with
$a^0\geq 0$ and $\ a^i>  0$ if $i>0$, and such that either $a^0\neq 0$ or
$a^i\neq 1$ for some $i\neq 0$, (i.e. those indices for which $M(n,\alpha)$
occurs in $V$) to sequences of length $n+2$, defined in the following way: If
$a^0\neq 0$ then $$
\sigma(\alpha)=(0,a^0,\ldots,a^n).
$$
 Otherwise, let $a^i$ be the first
element in the sequence, different from 1. Then
$$
\sigma(\alpha):=(a^0,\ldots,a^{i-1},1,a^i-1,a^{i+1},\ldots,a^n)=
(0,1,\ldots,1,1,a^i-1,a^{i+1},\ldots,a^n).
$$

Let
$$
M(n,\alpha):=\Gamma(Y\times
X(n),\L^{\otimes a_0}\times
\ldots \times \L^{\otimes a_n})$$ and recall that the partial differential
$d_n^{ii+1}$is a map
$$M(n,a^0,\ldots,a^n)\to M(n-1,a^0,\ldots,a^i+a^{i+1},\ldots,a^n).$$
Call temporarily the sequence
$(,a^0,\ldots,a^i+a^{i+1},\ldots,a^n)$  a {\it partial
differential}  of the sequence $(a^0,\ldots,a^n).$ The following easy lemma is
crucial in the construction of the homotopy.

\proclaim{Lemma 4.1} Suppose that
$\beta=(b^0,b^1,\ldots,b^{n-1})=(a^0,\ldots,a^i+a^{i+1},\ldots,a^n)$
is a partial differential of $\alpha=(a^0,\ldots,a^n).$
Then  $\sigma(\beta)$ is
a partial differential of $\sigma(\alpha)$. Also  $\alpha $ itself is a partial
differential of $\sigma(\alpha)$.    \endproclaim

\demo {Proof } The last statement is immediate from the definition of $\sigma$.
The remaining proof falls into several cases. If $a^0\neq 0$ then
clearly $b^0\neq 0$ and
$\sigma(\beta)=(0,b^0,b^1,\ldots,b^{n-1})$ is a partial
differential of   $\sigma(\alpha)=(0,a^0,\ldots,a^n).$ If
$a^0= 0$, but $b^0\neq
0$ then $i=0$ and
$\sigma(\beta)=(0,b^0,b^1,\ldots,b^{n-1})=\alpha$, which is a
partial differential of $\sigma(\alpha)$, by the last part of the lemma. If
finally $a^0=b^0=0$, let $a^l$ be the first non-zero element in the sequence
different from 1. There are now three cases. If $i<l$ then $1=a^i $ and
$a^i+a^{i+1}-1=a^{i+1}$ and hence $\sigma(\beta)=\alpha$. If $i=l$ then
$\sigma(\beta)=(0,1,\ldots,1,a^l+a^{l+1}-1,\ldots,a^n)$ is clearly a partial
differential of $\sigma(\alpha)=(0,1,\ldots,1,a^l-1,a^{l+1},\ldots,a^n).$
If $i>l$
then $\sigma(\beta)=(0,1,\ldots,1,a^l-1,\ldots,a^i+a^{i+1},\ldots,
a^n)$ and this is a partial differential of
$\sigma(\alpha)=(0,1,\ldots,1,a^l-1,
\ldots,a^n).$\enddemo

We will now construct a degree +1 map
$s:V.\to N.$, such that the inclusion map $V.\to N.$ is
homotopic to 0 through s,
i.e.  $v=dsv+sdv,\ v\in V$. This map will be constructed
inductively. Assume the
condition on models of the theorem . Note that a model complex $K(n,\L^{\otimes
a_0}\times\ldots\times \L^{\otimes a_{n+1}})$ is a subcomplex of $K.$, and
that it contains in the middle degree precisely those $M(n,\beta)$ for which
$\beta $ is a partial differential of $\alpha$.

\subheading {First step} To construct $$s_0 :V_0=\bigoplus_{i\geq
0}\Gamma(Y,\L^i)\to K.,$$ consider the model complex $K:=K(0,\L\times \L^{i-1})
$  (note that $i\geq 2$) $$ K:\ \ \Gamma(Y\times X,L\times \L^{i-1})\overset
d_1{^01}\to\to\Gamma(Y,L^i)\overset 0\to\to 0$$
By assumption this complex is acyclic in degree 0, hence $d_1^{01}$
is surjective and there is a splitting $s_0$. Choose
such a splitting $s_0$ for each $i\geq 2$.  Clearly $d_{12}s_0a=a$, and hence
$s_0$ is acceptable as the first part of a homotopy of the
inclusion of $V\subset
K.$ with the zero map.

\subheading {The induction step} Assume that $s_{k},\ k< n$ has already been
constructed, such that $s_k(M(k,\gamma))\subset M(k+1,\sigma(\gamma))$. If
$m\in K_n$, the relation $$d_nm=s_{n-2}d_{n-1}d_nm+d_ns_{n-1}d_nm =
d_ns_{n-1}d_nm$$ is valid. Hence in particular
$d_n(m-s_{n-1}d_nm)=d_nm-d_ns_{n-1}d_nm=0$, so
$m-s_{n-1}d_nm$ is a cycle. To define $s_n$, it is enough to consider
$ m\in M(n,\alpha)\subset V_n$. Consider the model complex
$K(n,\sigma(\alpha))\subset K.$ By the induction hypothesis and Lemma 4.1 it is
clear that $s_{n-1}d_nm$ is in $K(n,\sigma(\alpha))$ and hence also the
cycle $m-s_{n-1}d_nm$.
By the acyclicity of $K(n,\sigma(\alpha))$ in degree $n$ there is a splitting
$\theta : Z_n(K(n,\sigma(\alpha)))\to
K(n,\sigma(\alpha)))_{n+1}=M(n,\sigma(\alpha))$ of $d_n$.
Define  $$s_nm:=\theta(m-s_{n-1}d_nm)\in M(n,\sigma(\alpha))$$ Clearly
$m=d_{n+1}s_nm+s_{n-1}d_nm$, so $s_n$ is the next step of
the homotopy as desired.

\subheading{5.Remaining proofs}

 We will now prove the results stated in the introduction.
To get Theorem 4, just recall that the arguments of Mehta-Ramanathan (see
e.g. [R,1.13]), show that restriction of ample invertible sheaves to compatibly
Frobenius split subvarieties induces
surjections on global sections, so that ii-iii) of Lemma
4.3.2 are satisfied, for
ample $\L_\alpha$, if the corresponding unions of partial
diagonals are compatibly
split in $X(n+1)$. By the same argument, if $Y\times X(n)$ is Frobenius split
then  i) of the lemma is true. Hence it suffices to note that by elementary
properties of ampleness $\L_\alpha$ is ample when $\L$ is ample. Thus Theorem
5 gives Theorem 4.
Now note that it is
possible to restrict oneself to positive characteristics as described e.g.in
[R,3.11] and [B,1.3] in order to show that the hypothesis of Theorem 5 is
fulfilled in arbitrary characteristics.

{}From Theorem 4, (taking $Y=\emptyset$ ) Theorem 1 follows directly for
all homogeneous coordinate rings belonging to ample line bundles of those toric
varieties, for which the requisite splitting properties of
partial diagonals are
known. It is proven in [B, 3.1] that all homogeneous coordinate rings of the
type in the theorem are algebra retracts of such rings and hence Theorem 1
follows in complete generality by Lemma 2.3 a).

 Theorem 2 on the Koszul property of homogeneous coordinate
rings of Schubert varieties clearly follows from the relative version given
in Theorem 3, by applying this theorem first to $Y=\emptyset$ to obtain that
$R=R(X,\L)$ is a Koszul algebra and then another time to get that the map
$\R\to S$ is 1-linear and finally using  Lemma 2.3 c) to see
that $S$ is a Koszul
algebra.

Finally Theorem 3 for  an {\it ample} line bundle $\L$ in
positive characteristic,
 is a consequence of Theorem 5  using the fact that all the
partial diagonals are
compatibly split. That this is true is a
straightforward generalization, which we omit, of the ingenious proof of this
result for $n=3$ contained in Theorem 3.5 iii) of [R].
The extension of this case of Theorem 3 to a merely effective line
bundle $\L$, is immediate from the fact, proven e.g. in
[R,3.8, 3.11], that there
is a parabolic subgroup $Q$ such that, denoting the map $G/P\to G/Q$ by $\pi$,
the following  is true. First $\L=\pi^*(\L')$ where $\L'$ is an ample
line bundle on $G/Q$, and secondly $R(G/P,\L)\cong
R(G/Q,\L')$ and $R(Y,\L)\cong
R(\pi(Y),\L')$.
\bigskip

%\font\smc=cmcsc10

\long\def\clause#1#2{\par\smallbreak\hangafter=1\hangindent45pt\noindent
{\hbox to45pt{{#1\hfill}}{#2}}\hfill\par
\ifdim \lastskip <\smallskipamount \removelastskip \penalty 55\smallskip \fi}
%standard {\"a}r 30pt eller kanske 25.

\def\rit#1 #2: #3;#4;#5 #6 #7 \par{\noindent\clause
{[#1]}{{\smc #2, \sl #3, \rm #4 \bf #5} (#6), #7.}\bigskip}
\def\vrit#1 #2: #3;#4 \par{\noindent\clause
{[#1]}{{\smc #2, \sl #3,} #4.}\bigskip}

%\vfill \eject

\centerline{{\smc References}} \frenchspacing\bigskip\medskip

\vrit Ba J.Backelin: Some homological properties of "high" Veronese
subrings;J.Algebra 146(1992)1-17

\vrit{BF}   J. Backelin, R. Fr\"oberg: Koszul algebras, Veronese
subrings and rings with linear resolutions;Rev. Roumaine Math. Pures
Appl. 30 \break
 (1985), 85-97

\vrit B R.B\"o gvad: On the homogeneous ideal of a projective nonsingular
toric variety; University of Stockholm preprint,1994 (to appear in J. Tohoku
Math.)

\vrit CE  H.Cartan and S.Eilenberg:  Homological algebra ;
Princ.UP,1960

\vrit EG D.Eisenbud and S.Goto: Linear free resolutions and
minimal multiplicity;
J.Algebra 88(1984),89-133

\vrit {K1}  G.R.Kempf:  Some wonderful rings in algebraic geometry;
J.Algebra, vol 134, 1990

\vrit {K2}  G.R.Kempf: Wonderful rings and awesome modules;
Free resolutions in
commutative algebra and algebraic geometry, Sundance 90; Jones and Bartlett
Publishers, Boston London 1992

\vrit KR  G.R.Kempf and A.Ramanathan: Multi-cones over Schubert
varieties;Inventiones Math 1987

\vrit Le J.-M. Lemaire : Alg\'ebres connexes et homologie des espaces de
lacets;(Lect. Notes Math., vol. 422), Berlin Heidelberg New York: Springer 1974

\vrit L G.Levin: Large homomorphisms of local rings;
Math.Scand 46(1980), 209-215

\vrit  M S.Maclane:  Homology;
Berlin Heidelberg New York: Springer, 1963

\vrit  {R}   A. Ramanathan: Equations defining Schubert
varieties and Frobenius splitting of diagonals;
Publ.Math., Inst.Hautes \acuteaccent Etud.Sci. 65 (1987), 61-90

\vrit {Ro} J.-E.Roos: Commutative non-Koszul algebras having a linear
resolution of arbitrarily high order.Applications to torsion in loop space
homology;C.R.Acad.Sci.Paris Se'r.I Math.316 (1993), 1123-1128

\bigskip

DEPARTMENT OF MATHEMATICS

UNIVERSITY OF STOCKHOLM

S-106 91 STOCKHOLM

SWEDEN

{\bf email: rikard\@matematik.su.se}

\bye